\documentclass[11pt]{amsart}

\usepackage{epsfig,amssymb,latexsym,amscd,xypic}

\newtheorem{defi}{Definition}

\newtheorem{theo}{Theorem}
\newtheorem{lema}{Lemma}
\begin{document}

\title{Radford's formula for biFrobenius algebras and applications} 

\author{Walter Ferrer Santos and Mariana Haim}
\thanks{The first author would like to thank, Csic-UDELAR and
  Conicyt-MEC.}

\begin{abstract}
In a biFrobenius algebra $H$, in particular in the case that $H$ is a
finite dimensional Hopf algebra, the antipode $\mathcal S:H
\rightarrow H$ can be decomposed as $\mathcal S= _t\!\!c {\scriptstyle\circ}\/ c_\phi$ where
$c_\phi:H \rightarrow H^*$ and $_tc:H^* \rightarrow H$ are the
Frobenius and coFrobenius isomorphisms. We use this
decomposition to present an easy proof of Radford's formula for
$\mathcal S^4$. Then, in the case that the map $\mathcal S$ satisfies the
additional condition that $\mathcal S \star\, \operatorname{id} = \operatorname
{id} \star\,\, \mathcal S = u \varepsilon$, we prove the {\em trace \!\!\!
  formula}\/: 
$\operatorname {tr}(\mathcal S^2)=\varepsilon(t)\phi(1)$ 
. We finish by applying the above results
to study the semisimplicity and cosemisimplicity of $H$. 
\end{abstract}
\maketitle
\section{Introduction}

The concept of biFrobenius algebra (or simply bF algebra) was introduced by
Y. Doi and M. Takeuchi in \cite{kn:dt} as a generalization of finite
dimensional Hopf algebras. In particular, bF algebras are equipped
with a special antimorphism of algebras and coalgebras that plays the
role of the antipode for Hopf algebras and is denoted by the letter
$\mathcal S$. The theory of bF algebras was further
developped in \cite{kn:doi1}, \cite{kn:doi2}, \cite{kn:FS} and
\cite{kn:H}.  For the general theory of Hopf algebras we refer the reader to: 
\cite{kn:montgomery}, \cite{kn:schneider} and \cite{kn:sweedler}.

The main purpose of this paper is to begin the study of a class of
bF algebras introduced in \cite{kn:H}, that we call here the
SbF algebras. It consists of bF algebras with the
additional condition that the map $\mathcal S$ is the convolution
inverse of the identity in $H$.

A crucial role in the study of the semisimplicity and cosemisimplicity
of finite dimensional Hopf algebras is played by Radford's formula for
$\mathcal S^4$. Hence, after briefly recalling the basic definitions of bF
and SbF algebras in Section 2, we present in Section 3 a very short
proof of Radford's formula for bF algebras. In Section 4 we prove
that the usual formula for $\operatorname {tr}(\mathcal S^2)$ is
valid in the SbF situation. In Section 5 we study the semisimplicity
and cosemisimplicity of SbF algebras over fields of characteristic
zero.

The following notations will be in force along the paper. If $V$ is a
$\Bbbk$--linear space, $f\in V^*$ and $v \in V$ we define the
linear transformation $f | v : V \rightarrow V$ as:
$(f|v)(w)=f(w)v$. If $V$ is finite dimensional, the map $f \otimes v
\mapsto f|v: V^* \otimes V \rightarrow \operatorname{End}(V)$ is an
isomorphism. Frequently, we represent the elements of
$\operatorname{End}(V)$ as elements of $V^* \otimes V$ using the above
isomorphism. In this case if $T= \sum f_i \otimes v_i \in
\operatorname{End}(V)$, then $\operatorname{tr}(T)=\sum f_i(v_i)$.

\section{BiFrobenius algebras}

We review here the basic results and definitions of biFrobenius
algebras following \cite{kn:dt}.

Assume that $(H,m,1)$ is an associative unital $\Bbbk$--algebra and
$(H,\Delta,\varepsilon)$ is a coassociative, counital $\Bbbk$--coalgebra. We
use Sweedler's notation for the comultiplication.

We consider $\rightharpoonup:H \otimes H^* \rightarrow H^*$ and
 $\leftharpoonup:H^* \otimes H \rightarrow H^*$ --the standard left and
right actions of $H$ on $H^*$-- defined as
$(x \rightharpoonup \lambda)(y)=\lambda(yx)$ and $(\lambda \leftharpoonup x)(y)=\lambda(xy)$  
for $x,y \in H$, $\lambda \in H^*$.

Dually, one can take $\rightharpoonup:H^* \otimes H \rightarrow H$  and
$\leftharpoonup:H \otimes H^* \rightarrow H$ --the left and right action
of $H^*$ on $H$-- defined as follows: if $x \in H$, $\lambda \in H^*$: $x
\leftharpoonup \lambda=\sum \lambda(x_1)x_2$ and $\lambda\rightharpoonup x = \sum
x_1\lambda(x_2)$.

The algebra $H$ is called a Frobenius algebra if there exists an
element $\phi \in H^*$ such that $H^*=\phi \leftharpoonup H$ (which
implies that $H^*=H \rightharpoonup \phi$).

The coalgebra $H$ is called a coFrobenius coalgebra if there exists an
element $t \in H$ such that $H=t\leftharpoonup H^*$ (which implies that 
$H=H^* \rightharpoonup t$).
 
If $H$ is Frobenius or coFrobenius, it is finite dimensional. 

If $H$ is Frobenius and coFrobenius, $1 \in H$ is a group like element and
$\varepsilon:H \rightarrow \Bbbk$ is an algebra homomorphism, $\phi$
and $t$ can be assumed to satisfy the additional conditions: $t
\leftharpoonup \phi =1$, $\phi \leftharpoonup t = \varepsilon$. In
that case $t$ is a right integral and $\phi$ a right
cointegral of $H$ respectively. In other words for all $x \in H$, $tx=\varepsilon(x)t$
and $\phi(x)1=\sum \phi(x_1)x_2$. Notice also that in this situation
$\phi(t)=1$.

In the above context, one can consider the map $\mathcal S:H
\rightarrow H$, $\mathcal S(x)= \sum \phi(t_1x)t_2$ that is called the
antipode morphism. In particular $S(1)=t\leftharpoonup \phi=1$ and
$\varepsilon {\scriptstyle\circ}\/ \mathcal S=\phi \leftharpoonup t = \varepsilon$. 

In this situation we call $_tc,c_t:H^* \rightarrow H$ the
$\Bbbk$--linear maps $_tc(f)=t \leftharpoonup f$ and $c_t(f)= f
\rightharpoonup t$. 

Similarly define  $_\phi c , c_\phi:H \rightarrow H^*$ as the
$\Bbbk$--linear maps $_{\phi}c(x)=\phi \leftharpoonup x$ and $c_{\phi}(x)= x
\rightharpoonup \phi$. 

In this notation $\mathcal S= _t\!\!c {\scriptstyle\circ}\/ c_{\phi}$ is a linear
bijective endomorphism of $H$. The composition inverse of $\mathcal S$
is denoted as $\overline{\mathcal S}$.

In \cite{kn:dt} --the notations and hypotheses are the ones considered above--
the authors define the 7--uple $(H,m,1,\Delta,\varepsilon,t,\phi)$ to
be a biFrobenius algebra if $\mathcal S$ is an antiautomorphism of 
algebras and coalgebras.  

In this case $\phi {\scriptstyle\circ}\/ \mathcal S =\lambda$ and ${\mathcal S}t=s$
are respectively a left cointegral and a left integral for $H$. Moreover, one
can define the elements $a \in H$ and $\alpha \in H^*$ 
--the modular element and the modular function-- as $\alpha = t
\rightharpoonup \phi$ and $a= \phi \rightharpoonup t$.

In this case, $a$ is a group like element, $\alpha$ a morphism of
algebras, ${\mathcal S}a=a^{-1}$, $\alpha {\scriptstyle\circ}\/ \mathcal S =
\alpha^{-1}$, and for all $x \in H$: $sx=\alpha^{-1}(x)s$ 
and $\lambda(x)a^{-1}=\sum \lambda(x_1)x_2$

Moreover --see \cite{kn:dt}-- the maps $x \mapsto ax\, , \, x
\mapsto xa: H \rightarrow H$ are automorphisms of coalgebras,
i.e., for all $x \in H$, 
\begin{equation} \sum (ax)_1 \otimes (ax)_2=\sum ax_1 \otimes ax_2
\end{equation}
and
\begin{equation} 
\sum (xa)_1 \otimes (xa)_2=\sum x_1a \otimes x_2a 
\end{equation}

Similarly, the maps $x \mapsto (x \leftharpoonup \alpha) \, , \, x
\mapsto (\alpha \leftharpoonup x): H \rightarrow H$ are morphisms of
algebras. In other words for all $x,y \in H$:
\begin{equation} \sum\alpha((xy)_1)(xy)_2 = \sum \alpha(x_1)\alpha(y_1)x_2y_2
\end{equation}
and
 \begin{equation} \sum (xy)_1\alpha((xy)_2) = \sum x_1y_1\alpha(x_2)\alpha(y_2)
\end{equation}

Using the uniqueness of the integrals one can prove that: $a
\rightharpoonup \phi = \lambda$\,, $\phi(\mathcal S t)=
\lambda(t)=\phi(s)=\alpha(a)=1$.

\noindent
The following formul\ae\/ are valid in an arbitrary
  bF-algebra.

\begin{enumerate}
\renewcommand{\theenumi}{\roman{enumi}}
\item For all $x \in H$, $x=\sum \phi(t_1x)\overline{\mathcal
  S}(t_2)$  
\item For all $x,y \in H$,  $\sum \phi(y_1 x)y_2= \sum \phi(y
  x_1)\mathcal S(x_2)$
\item For all $x,y \in H$, $\sum \phi(xy_2)y_1= \sum \phi(x_2
  y)\mathcal S(x_1)a$
\item If we put in (iii) $x=t$ and $y=t$ we have: 
\begin{enumerate} 
\item $\mathcal S(\alpha \rightharpoonup x)a =\sum \phi(x t_2)t_1$
\item $a \overline {{\mathcal S}}(x) =\sum \phi(t_2x)t_1$
\end{enumerate}
\item If we write in equation (iii) $x=\mathcal S(t)$ we have: \\ $\sum
  \phi(\mathcal St_1x) t_2 = \overline{\mathcal S}^2(\alpha^{-1}
  \rightharpoonup x)a^{-1}$
\end{enumerate}

For example equation (ii) can be deduced from the fact that $\mathcal
S$ is an antimorphism of coalgebras as follows. 

Form the equality $y=\sum \phi(y {\overline {\mathcal S}}t_2)t_1$, we
deduce that 
\begin{equation} \label{eqn:primera}\sum y_1x \otimes y_2= \sum \phi(y {\overline
    {\mathcal S}}t_3)t_1x \otimes t_2
\end{equation}

As $\mathcal S(x)= \sum \phi(t_1x)t_2$, we deduce that $\mathcal
S(x_2) \otimes \mathcal S(x_1)=
\sum \phi(t_1x)t_2 \otimes t_3$ and

\begin{equation} \label{eqn:segunda} \sum yx_1 \otimes {\mathcal S}x_2 = \sum (y {\overline{\mathcal
      S}}t_3 \otimes t_2)\phi(t_1x)
\end{equation}
 
Applying $\phi \otimes \operatorname {id}$ to the equations
 (\ref{eqn:primera}) and (\ref{eqn:segunda}) we get (ii). 

For later use we compute $c_{\phi}{\scriptstyle\circ}\/ \,_tc: H^* \rightarrow H^*$. 

If $\gamma \in H^*$ and $x \in H$ we have

\begin{equation} \label{eqn:tercera} (c_{\phi} {\scriptstyle\circ}\/ \,_tc)(\gamma)(x)=(_tc(\gamma)
  \rightharpoonup \phi)(x)=\phi(x (_tc(\gamma)))=
\end{equation} 
\begin{equation*} \sum \gamma(t_1)\phi(xt_2)=\gamma(\sum
  \phi(xt_2)t_1)=\gamma(\mathcal S(\alpha \rightharpoonup x)a).
\end{equation*}

The last equality follows from (iv).

For a bF algebra, it is useful to consider the Nakayama and coNakayama
automorphisms --see \cite{kn:dt} for the situation of bF algebras and
\cite{kn:schneider} for the case of finite dimensional Hopf algebras--. This morphisms are denoted as
$\mathcal N , \, ^c\mathcal N: H \rightarrow H$ and are defined by the
following equalities:
\begin{equation} \phi \leftharpoonup x = \mathcal N(x) \rightharpoonup
  \phi, \,\,\, \text{for all} \,\,\, x \in H
\end{equation}

\begin{equation} t \leftharpoonup f = (f {\scriptstyle\circ}\/\/ ^c\/\mathcal N)
  \rightharpoonup t, \,\,\, \text{for all} \,\,\, f \in H^*
\end{equation}

In more explicit terms $\mathcal N$ and $^c\mathcal N$ can be
characterized by the following equations: for all $x,y \in H$,
$\phi(xy)=\phi(y \mathcal N(x))$ and $\sum\/ ^c\!\mathcal N(t_2) \otimes
t_1= \sum t_1 \otimes t_2$ if $\Delta(t)=\sum t_1 \otimes t_2$. 

It is easy to show that $\mathcal N$ is an automorphism of algebras and
that $^c\mathcal N$ is an automorphism of coalgebras.

If we apply $^c\mathcal N^{-1}$ to  the equality (iv) we obtain that 
$^c\mathcal N^{-1}(a \overline {\mathcal S}(x))=\sum \phi(t_2x)\/\,\,
^c\mathcal N^{-1}(t_1)=\sum \phi(t_1x)t_2$. 

Hence, we conclude that:

\begin{equation} \label{eqn:conak} ^c\mathcal N(x)=a\overline{\mathcal S}^2(x)
\end{equation}

Similarly using (iv) again, we deduce that $a (\overline {\mathcal
  S}\mathcal N)(x)= \sum \phi(t_2 \mathcal N(x))t_1= \sum
  \phi(xt_2)t_1= \mathcal S(\alpha \rightharpoonup x)a$.

Hence, we conclude that:

\begin{equation} \label{eqn:nak} \mathcal N(x)= a^{-1}\mathcal S^2(\alpha
  \rightharpoonup x)a
\end{equation}

In particular $^c\mathcal N(a)=a^2$ and $\mathcal N(a)=a$. 

Next we compute the traces of $\mathcal N$ and $^c\mathcal N$.

From the equality $\mathcal S(x)=\sum \phi(t_1x)t_2$ we obtain that
$^c\mathcal N \mathcal S x=\sum \phi(t_2x)t_1$ and $\overline
{\mathcal S }\/\,^c\!\mathcal N \mathcal S=\sum \phi \leftharpoonup t_2
\otimes \overline {\mathcal S}t_1$.

Then 
\begin{equation} \label{eqn:trazaconak} \operatorname {tr} (^c\mathcal N)= \operatorname {tr}
  (\overline {\mathcal S}\/\,^c\!\mathcal N \mathcal S)= \sum \phi(t_2 \overline
  {\mathcal S}t_1)=(\phi {\scriptstyle\circ}\/ \overline {\mathcal S})((\operatorname
  {id}\star\/\,\, \mathcal S)(t))
\end{equation}  

Similarly from $x=\sum \phi(t_1x)\overline {\mathcal S}t_2$ we deduce
that $\mathcal N x = \sum \phi(xt_1) \overline{\mathcal S}t_2$, and
then:

\begin{equation}\label{eqn:trazanak} \operatorname {tr} (\mathcal N)=
  \phi((\overline {\mathcal S}t_2)t_1)=(\phi {\scriptstyle\circ}\/ \overline {\mathcal
  S})((\mathcal S \star \operatorname {id})(t))
\end{equation}

We need one more trace computation.

The equality  $x=\sum \phi(t_1x)\overline {\mathcal S}t_2$ can be
written as: $\operatorname {id} = \sum \phi \leftharpoonup t_1 \otimes
\overline {\mathcal S}t_2$.

\begin{equation}\label{eqn:trazaid} \operatorname {dim} (H)=
  \operatorname {tr}(\operatorname {id})=\sum \phi(t_1 \overline {\mathcal
  S}t_2) = \phi((\operatorname {id} \star\/ \,\,\overline {\mathcal S})(t))
\end{equation}

\medskip
\medskip

Next --following \cite{kn:H}-- we define a special family of bF
algebras, whose representation theoretical properties can be put under
a stricter control than for the general bF algebras.
In \cite{kn:H}, the second author of this paper
considers the class of biFrobenius algebras
satisfying the additional condition that 
\begin{equation} \label{eqn:inviden} \mathcal S \star \,
\operatorname {id}=\operatorname {id} \star \,\,\mathcal S =
u\varepsilon.
\end{equation} 

Clearly not all bF algebras satisfy condition
(\ref{eqn:inviden}) --see for example \cite{kn:dt}-- and in the 
mentioned paper \cite{kn:H}  
using known results on the existence of large Hadamard
matrices, a family of bF algebras of arbitrarily large
dimension satisfying condition (\ref{eqn:inviden}) and that are not Hopf
algebras is constructed. 

In other words, for the objects of this family of bF
algebras, the multiplication and comultiplication are not related by
the so called pentagonal axiom but only by the weaker condition
(\ref{eqn:inviden}). 

It is convenient to give an explicit name to the family of all bF
algebras satisfying the condition (\ref{eqn:inviden}).

\begin{defi} If $H$ is a biFrobenius algebra, we say that $H$ is of
  type $S$ --or that $H$ is an SbF algebra--, if the map $\mathcal
  S$ is the convolution inverse of the identity in $\operatorname
  {End}_{\Bbbk}(H)$.
\end{defi}

{\large {\bf Observation 1.}}\/\,\, It is worth noticing that one can
construct bF algebras $H$ with the property that the identity map 
$\operatorname {id} : H \rightarrow H$ is convolution invertible, 
but that are not SbF algebras. In this case the
convolution inverse of the identity is not the antipode $\mathcal
S=\sum \phi(t_1x)t_2$ of the bF algebra. 

An example of the above situation is the following. Let $\Bbbk$ be a
field of characteristic different from 2 and $H$ the vector space 
linearly generated by three elements $H= \langle 1,x,y \rangle$. We
endow $H$ with the following structure of bF algebra. 
The multiplication table of $H$ is characterized by the conditions: 
$1 \in H$ is the unit element and 

$$
x^2=\frac{1}{2}x + \frac{3}{2} y, \ \ \ xy=yx=2+\frac{1}{2}x +
\frac{1}{2}y, \ \ \ y^2=\frac{3}{2}x+\frac{1}{2}y.
$$ 

The comultiplication is given as: $$ \Delta(x) = \frac{1}{2} x \otimes
x, \\\ \Delta(y) = \frac{1}{2} y \otimes y.$$  Moreover $1 \in H$ is a
group like element and $\varepsilon$ is defined as:
$\varepsilon(x)=\varepsilon(y)=2$.
  
If we take  $t=1+x+y$ and $\phi(a +bx+cy)=a$ one can show directly
that $(H,m,1,\Delta,\varepsilon,t,\phi)$ is a bF algebra with antipode
$\mathcal S : H \rightarrow H$ given by: $\mathcal S(1)=1, \mathcal
S(x)=y, \mathcal S(y)=x$. 

Moreover a direct computation shows that the map $\Sigma: H
\rightarrow H$ given as $$\Sigma(1)=1, \\\ \Sigma (x)= -\frac{2}{3} -
\frac{2}{3}x + 2y, \\\ \Sigma (y) =  -\frac{2}{3} + 2x -
\frac{2}{3}y$$ is the convolution inverse of the identity in $H$.   

Notice that the above example is in fact a group like algebra in accordance
with the definition of \cite{kn:doi2}

\section{A short proof of Radford's formula}

Radford's formula for $\mathcal S^4$ in the case of Hopf algebra 
was first proved in full
generality in \cite{kn:r} and with predecessors in \cite{kn:s} and
\cite{kn:l}. A more recent proof, that is in  the spirit of the one 
we present below, appears in \cite{kn:schneider}. 
Generalizations of the formula
from the case of Hopf algebras to other situations, braided Hopf algebras,  bF
algebras --braided and classical, quasi Hopf algebras, weak Hopf
algebras, Hopf algebras over rings, and even for the very general case
of finite tensor categories, 
can be found in the following references: \cite{kn:bkl}, \cite{kn:dt},
\cite{kn:eno},\cite{kn:hn}, \cite{kn:ks}, \cite{kn:k} and
\cite{kn:n}. 

\medskip

For the proof we write $\mathcal S^2=\/_tc {\scriptstyle\circ}\/ (c_{\phi}{\scriptstyle\circ}\/\/ _tc){\scriptstyle\circ}\/
c_{\phi}$. In other words and more explicitly, we have that:
\begin{equation}
\mathcal S^2(x)=\sum ((c_{\phi}{\scriptstyle\circ}\/\/_tc)(x \rightharpoonup
\phi))(t_1)t_2.
\end{equation}

Then, from equation (\ref{eqn:tercera}) we deduce that

\[\mathcal S^2(x)=\sum (x \rightharpoonup \phi)(\mathcal S(\alpha
\rightharpoonup t_1)a)t_2= \sum \phi(({\mathcal
  S}t_1)ax)\alpha(t_2)t_3.\]

Moreover:  

\begin{equation} \sum \phi(({\mathcal S}t_1)ax)\alpha(t_2)t_3=\sum
\alpha^{-1}(\phi(({\mathcal S}t_1)ax)\mathcal S(t_2))t_3=
\end{equation}
\begin{equation*}\sum
\alpha^{-1}(\phi({\mathcal S}(t_1)_2ax){\mathcal S}(t_1)_1)t_2 =
\sum \alpha^{-1}(a \phi({\mathcal S}t_1 a x_2)\overline{\mathcal S}(ax_1))t_2
\end{equation*}

where in the last equality we used (iii).

Thus, $\mathcal S^2(x)= \sum \phi({\mathcal S}t_1 a
x_2)\alpha(x_1)t_2= \sum \alpha(x_1) {\overline{\mathcal
    S}}^2(\alpha^{-1} \rightharpoonup (ax_2))a^{-1}$, where this last
equality follows from (v).

Moreover, $$\mathcal S^2(x)={\overline{\mathcal S}}^2(\sum
\alpha(x_1)ax_2 \alpha^{-1}(x_3)a^{-1})={\overline{\mathcal
    S}}^2(a(\alpha^{-1}\rightharpoonup x \leftharpoonup
\alpha)a^{-1})$$ and then
\[ {\mathcal S}^4(x)= a(\alpha^{-1}\rightharpoonup x \leftharpoonup
\alpha)a^{-1}\]

In particular if $\alpha = \varepsilon$ and $a=1$, i.e., if $H$ is
unimodular and counimodular, the antipode satisfies the equality
${\mathcal S}^4= \operatorname {id}$. 

\section{The trace formula}

From the equality $\mathcal S(x)=\sum \phi(t_1x)t_2$ we obtain that:

\begin{equation} {\mathcal S}^2(x)=\sum \phi(t_1x)\mathcal S t_2 
\end{equation}
and
\begin{equation}
{\mathcal S}^2= \sum \phi \leftharpoonup t_1 \otimes
  \mathcal S t_2 .
\end{equation}

\begin{theo} In the situation above, if $H$ is a bF algebra, then
  $\operatorname {tr}({\mathcal S}^2)=\phi((\mathcal S \star\/
  \operatorname {id})(t))$. Moreover, if $H$ is of type S, then
  $\operatorname {tr}({\mathcal S}^2)=\varepsilon(t)\phi(1)$. 
\end{theo}
\proof Taking traces in the equality ${\mathcal S}^2= \sum \phi
\leftharpoonup t_1 \otimes \mathcal S t_2$, we deduce that 
$\operatorname {tr}({\mathcal S}^2)=\sum (\phi
\leftharpoonup t_1)(S t_2)=\phi((\mathcal S \star\/
  \operatorname {id})(t))$
\qed

\medskip 
\medskip

{\large {\bf Observation 2.}}\/\,\,In \cite{kn:dt}, the authors present the
  example that they call ${\operatorname B}_4$, of a bF algebra
  defined as follows. As an algebra ${\operatorname
  B}_4=\Bbbk[X]/(X^4)$. If we call $x = X + (X^4)$ and consider the
  basis $\mathcal B=\{1,x,x^2,x^3\}$, then the coalgebra structure of
  ${\operatorname B}_4$ is given by the 
following rules: 1 is a group like element, $x$ and $x^2$  are
  primitive elements, and $\Delta(x^3)= 1 \otimes x^3 + x \otimes x^2
  + x^2 \otimes x + x^3 \otimes 1$. Moreover $\varepsilon$ is given as
  $\varepsilon(1)=1,
  \varepsilon(x)=\varepsilon(x^2)=\varepsilon(x^3)=0$.

If we call ${\mathcal B}^*=\{1^*,x^*,x^{2*},x^{3*}\}$ the dual basis
of $V^*$, it is easy to show that $t=x^3$ and $\phi=x^{3*}$.

In this situation $\mathcal S= \operatorname{id}$,
$\phi(1)=\varepsilon(t)=0$ and $\operatorname {tr}(\mathcal
S^2)=\operatorname{tr}(\operatorname {id})=4$. 

This example shows that the $S$--condition is crucial for the validity
of the trace formula. Observe also that $\operatorname {B}_4$ is not
semisimple or cosemisimple, but $\mathcal S^2=\operatorname {id}$, see
\cite{kn:H}. For Hopf algebras a classical theorem due to Larson and
Radford (see \cite{kn:lr1} and \cite{kn:lr2}) guarantees that in
the case of characteristic zero, if $H$ is semisimple and
cosemisimple, then $\mathcal S^2= \operatorname {id}$. 

\section{Semisimplicity and cosemisimplicity of bf algebras.}

The methods we present in this section are similar to the ones
appearing in \cite{kn:m} and \cite{kn:schneider}.

In the case that $H$ is unimodular we deduce from equation
(\ref{eqn:nak}) that $\mathcal N(x)=a^{-1} \mathcal S^2(x) a$ and
hence the map $\mathcal N$ is also a morphism of coalgebras. 

If $H$ is counimodular, i.e. if $a=1$, then $^c\mathcal N = \overline
{\mathcal S}^2$ and in this case $^c\mathcal N$ is also a morphism of
algebras. 

In the case that $H$ is simultaneously unimodular and counimodular,
from Radford's formula we deduce that $\mathcal S^4= \operatorname
{id}$ and from the above considerations that: $\mathcal N =\/\,^c\mathcal
N=\mathcal S^2=\overline {\mathcal S}^2$ and $\mathcal
N^2=\operatorname {id}$.  

In this situation we also have that: ${\mathcal S}t=t$ and $\phi
{\scriptstyle\circ}\/
\mathcal S=\phi$.

Also: $\sum\mathcal S t_2 \otimes \mathcal S t_1 = \sum t_1 \otimes
t_2$, $\sum t_2 \otimes \mathcal S t_1 = \sum {\overline {\mathcal
    S}}t_1 \otimes t_2$ , $\sum \mathcal N t_2 \otimes t_1 = \sum t_1 \otimes
t_2$ and $\sum \mathcal S t_2 \otimes t_1 = \sum {\overline {\mathcal
    S}}t_1 \otimes t_2$. 

Then:

\begin{equation}\label{eqn:todot}
(\operatorname {id} \star\/\,\, \overline {\mathcal S})(t)= (\overline
  {\mathcal S} \star {\operatorname {id}})(t)= 
\end{equation}
\begin{equation*}
\sum ({\mathcal S} t_2)t_1=
  \sum t_2 ({\mathcal S} t_1)= \sum ({\overline {\mathcal S}} t_1)t_2 =
  \sum t_1 ({\overline {\mathcal S }}t_2)
\end{equation*}  

\begin{lema} \label{lema:main} If $H$ is a unimodular and counimodular biFrobenius
  algebra of type $S$, then $\operatorname {tr}(\mathcal
  N)=\phi(1) \varepsilon (t)$. Moreover, in the case that the base
  field has characteristic zero and 
$(\operatorname {id} \star\/\,\, \overline {\mathcal S})(t)=
  \varepsilon(t)1$, then $\mathcal N = \mathcal S^2 =\/ ^c\!\mathcal N =
  \operatorname {id}$. 
\end{lema}
\proof We already observed that in this situation $\mathcal N^2 =
\operatorname {id}$. From the equation (\ref{eqn:trazaid}) we deduce
that $\operatorname
{dim}(H)=\phi(\varepsilon(t)1)=\varepsilon(t)\phi(1)=\operatorname
{tr}(\mathcal N)$. Hence, as all the eigenvalues of $\mathcal N$ are $\pm 1$
and its sum --the trace of $\mathcal N$-- equals the dimension of $H$, we
conclude that the eigenvalue $-1$ cannot appear so that $\mathcal N=
\operatorname {id}$. 
\qed

If $H$ is a semisimple biFrobenius algebra, it is easy to show that
$\varepsilon(t) \neq 0$ --see for example \cite{kn:H}--. Applying
$\varepsilon$ to the equality $xt=\alpha(x)t$ that is valid for all $x
\in H$, we deduce that
$\varepsilon(x)\varepsilon(t)=\alpha(x)\varepsilon(t)$, and then that
$\alpha=\varepsilon$. In other words a semisimple bF algebra is
counimodular. 

Similarly if $H$ is cosemisimple, one can conclude that $H$ is
unimodular, i.e., $a=1$. 

The above results can be summarized in the following theorem.

\begin{theo} \label{theo:main} Assume that $H$ is a biFrobenius algebra of type $S$
  defined over an algebraically closed field of characteristic
  zero. If $\mathcal S^2 = \operatorname {id}$, then $H$ is semisimple
  and cosemisimple. Conversely, if $H$ is semisimple and cosemisimple,
  and $\sum (\mathcal S t_2)t_1 = \varepsilon(t) 1$, then $\mathcal S^2
  = \operatorname {id}$.
\end{theo}
\proof If $\mathcal S^2 =\operatorname {id}$, then $\operatorname
       {tr}(\mathcal S^2)=\phi(1)\varepsilon(t)=\operatorname
       {tr}(\operatorname {id})= \operatorname {dim}(H)$. In this
       situation, $\phi(1) \neq 0$ and $\varepsilon (t) \neq 0$. It
       is known that in this case --see for example \cite {kn:doi1} or 
\cite{kn:H}-- $H$ is semisimple and cosemisimple. The rest of the
       results follow from Lemma \ref{lema:main}.
\qed

\medskip

\noindent {\large {\bf Observation 3.}}\/\,\,It would be interesting to know if
--similarly than for
the situation of finite dimensional Hopf algebras-- the result
       is true without assuming the hypothesis that  $\sum (\mathcal S
       t_2)t_1 = \varepsilon(t) 1$.

\vspace*{1cm}
\noindent {\sc Walter Ferrer Santos}\\
Facultad de Ciencias\\
Universidad de la Rep\'ublica\\
Igu\'a 4225\\
11400 Montevideo\\
Uruguay\\
e-mail: {\tt wrferrer@cmat.edu.uy}

\vspace*{1cm}
\noindent{\sc Mariana Haim}\\
Facultad de Ciencias\\
Universidad de la Rep\'ublica\\
Igu\'a 4225\\
11400 Montevideo\\
Uruguay\\
e-mail: {\tt negra@cmat.edu.uy}

\end{document}